\documentstyle[11pt]{article}
\makeatletter
\ifx\@auxdefsloaded\relax\endinput
\else\let\@auxdefsloaded\relax\fi

\def\newenvironment{%
   \@ifnextchar *{\@@newenv{\global\@ignoretrue}}{\@@newenv{}*}}

\def\@@newenv#1*#2{%
   \@ifnextchar [{\@newenv{#1}{#2}}{\@newenv{#1}{#2}[0]}}

\long\def\@newenv#1#2[#3]#4#5{%
   \expandafter\newcommand\csname#2\endcsname[#3]{#4}%
   \expandafter\long\expandafter\def\csname end#2\endcsname{#5#1}}

\def\renewenvironment{%
   \@ifnextchar *{\@@renewenv{\global\@ignoretrue}}{\@@renewenv{}*}}

\def\@@renewenv#1*#2{%
   \@ifnextchar [{\@renewenv{#1}{#2}}{\@renewenv{#1}{#2}[0]}}

\long\def\@renewenv#1#2[#3]#4#5{%
   \expandafter\renewcommand\csname#2\endcsname[#3]{#4}%
   \expandafter\long\expandafter\def\csname end#2\endcsname{#5#1}}

\def\newoptcommand#1#2{%
   \@ifnextchar [{\@optargdef#1#2}{\@optargdef#1#2[1]}}

\def\renewoptcommand#1#2{%
   \edef\@tempa{\expandafter\@cdr\string#1\@nil}%
   \@ifundefined{\@tempa}{%
      \@latexerr{\string#1\space undefined}\@ehc}{}%
   \@ifnextchar [{\@reoptargdef#1#2}{\@reoptargdef#1#2[1]}}

\long\def\@optargdef#1#2[#3]#4{%
   \@ifdefinable #1{\@reoptargdef#1#2[#3]{#4}}}

\long\def\@reoptargdef#1#2[#3]#4{%
   \@tempcnta#3\relax \@tempcntb \@ne
   \let#1\relax \let\@tempa\relax
   \edef\@tempb{\long\def\csname\string#1\endcsname
      [\@tempa\the\@tempcntb]}%
   \advance\@tempcntb \@ne \advance\@tempcnta \m@ne
   \@whilenum\@tempcnta>0\do{%
      \edef\@tempb{\@tempb\@tempa\the\@tempcntb}%
      \advance\@tempcntb \@ne \advance\@tempcnta \m@ne}%
   \let\@tempa=##\@tempb{#4}%
   \def#1{\@ifnextchar [{\csname\string#1\endcsname}{%
      \csname\string#1\endcsname[#2]}}}

\def\newoptenvironment{%
   \@ifnextchar *{\@@newoptenv{\global\@ignoretrue}}{%
      \@@newoptenv{}*}}

\def\@@newoptenv#1*#2#3{%
   \@ifnextchar [{\@newoptenv{#1}{#2}{#3}}{%
      \@newoptenv{#1}{#2}{#3}[0]}}

\long\def\@newoptenv#1#2#3[#4]#5#6{%
   \expandafter\newoptcommand\csname#2\endcsname{#3}[#4]{#5}%
   \expandafter\long\expandafter\def\csname end#2\endcsname{#6#1}}

\def\renewoptenvironment{%
   \@ifnextchar *{\@@renewoptenv{\global\@ignoretrue}}{%
      \@@renewoptenv{}*}}

\def\@@renewoptenv#1*#2#3{%
   \@ifnextchar [{\@renewoptenv{#1}{#2}{#3}}{%
      \@renewoptenv{#1}{#2}{#3}[0]}}

\long\def\@renewoptenv#1#2#3[#4]#5#6{%
   \expandafter\renewoptcommand\csname#2\endcsname{#3}[#4]{#5}%
   \expandafter\long\expandafter\def\csname end#2\endcsname{#6#1}}

\newcounter{keepoptional}
\newcounter{optnestctr}

\newoptenvironment*{optional}{1}[1]{%
   \ifnum#1>\value{keepoptional}\relax
      \setcounter{optnestctr}{0}\@powerup\expandafter\@endopt\fi
   \ignorespaces}{}

\def\@powerup{\catcode`\{=12 \catcode`\}=12 \catcode`\\=12 \relax}
\def\@powerdown{\catcode`\{=1 \catcode`\}=2 \catcode`\\=0 \relax}
\begingroup \catcode`|=0 \catcode`[=1 \catcode`]=2 \@powerup
   |long|gdef|@endopt#1\end{optional}[%
      |@beginopt#1\begin{optional}*]%
   |long|gdef|@beginopt#1\begin{optional}[%
      |@ifstar[%
         |ifnum|value[optnestctr]>0|relax
            |addtocounter[optnestctr][-1]|let|@zzzap=|@endopt
         |else
            |@powerdown|end[optional]|let|@zzzap=|relax|fi
	 |@zzzap]%
      [
         |addtocounter[optnestctr][1]|@beginopt]]%
|endgroup

\makeatother
\makeatletter
\ifx\@auxdefsloaded\relax
\else \input{auxdefs.sty}\fi
\newskip\dgARROWLENGTH  \dgARROWLENGTH=2.5em\relax
\newskip\dgHORIZPAD     \dgHORIZPAD=1em\relax
\newskip\dgVERTPAD      \dgVERTPAD=2ex\relax
\newskip\dgLABELOFFSET  \dgLABELOFFSET=.7ex\relax
\newcount\dgARROWPARTS  \dgARROWPARTS=4\relax
\newcommand{\dgeverynode}{\displaystyle}
\newcommand{\dgeverylabel}{\scriptstyle}
\newskip\dgDOTSPACING   \dgDOTSPACING=0.35em
\newskip\dgDOTSIZE      \dgDOTSIZE=1.5\fontdimen8\tenln
\newcount\dgMAXSQUARE   \dgMAXSQUARE=4\relax
\newcount\dgMINDBLSQ    \dgMINDBLSQ=10\relax
\newskip\dgCOLUMNWIDTH  \dgCOLUMNWIDTH=2em\relax
\chardef\f@ur=4
\@namedef{dgo@}{\let\dg@VECTOR=\vector
   \dg@LBLPOS=\dgARROWPARTS \divide\dg@LBLPOS\tw@}%
\@namedef{dgo@-}{\let\dg@VECTOR=\line}%
\@namedef{dgo@!}{\let\dg@VECTOR=\dg@novector}%
\@namedef{dgo@1}{\dg@LBLPOS=\@ne\relax}%
\@namedef{dgo@2}{\dg@LBLPOS=\tw@\relax}%
\@namedef{dgo@3}{\dg@LBLPOS=\thr@@\relax}%
\@namedef{dgo@4}{\dg@LBLPOS=\f@ur\relax}%
\@namedef{dgo@5}{\dg@LBLPOS=5\relax}%
\@namedef{dgo@6}{\dg@LBLPOS=6\relax}%
\@namedef{dgo@7}{\dg@LBLPOS=7\relax}%
\@namedef{dgo@8}{\dg@LBLPOS=8\relax}%
\@namedef{dgo@9}{\dg@LBLPOS=9\relax}%
\def\dgt@e{\dg@DX=\@ne \dg@DY=\z@ \dg@SIZE=\@ne}%
\def\dgt@w{\dg@DX=\m@ne \dg@DY=\z@ \dg@SIZE=\@ne}%
\def\dgt@n{\dg@DX=\z@ \dg@DY=\@ne \dg@SIZE=\@ne}%
\def\dgt@s{\dg@DX=\z@ \dg@DY=\m@ne \dg@SIZE=\@ne}%
\def\dgt@ne{\dg@DX=\@ne \dg@DY=\@ne \dg@SIZE=\@ne}%
\def\dgt@se{\dg@DX=\@ne \dg@DY=\m@ne \dg@SIZE=\@ne}%
\def\dgt@nw{\dg@DX=\m@ne \dg@DY=\@ne \dg@SIZE=\@ne}%
\def\dgt@sw{\dg@DX=\m@ne \dg@DY=\m@ne \dg@SIZE=\@ne}%
\def\dgt@nne{\dg@DX=\@ne \dg@DY=\tw@ \dg@SIZE=\@ne}%
\def\dgt@nnw{\dg@DX=\m@ne \dg@DY=\tw@ \dg@SIZE=\@ne}%
\def\dgt@sse{\dg@DX=\@ne \dg@DY=-\tw@ \dg@SIZE=\@ne}%
\def\dgt@ssw{\dg@DX=\m@ne \dg@DY=-\tw@ \dg@SIZE=\@ne}%
\def\dgt@ene{\dg@DX=\tw@ \dg@DY=\@ne \dg@SIZE=\tw@}%
\def\dgt@ese{\dg@DX=\tw@ \dg@DY=\m@ne \dg@SIZE=\tw@}%
\def\dgt@wnw{\dg@DX=-\tw@ \dg@DY=\@ne \dg@SIZE=\tw@}%
\def\dgt@wsw{\dg@DX=-\tw@ \dg@DY=\m@ne \dg@SIZE=\tw@}%
\def\dggeometry{
   \dg@ZTEMP=\dg@XGRID \multiply\dg@ZTEMP\tw@
   \ifnum\dg@YGRID=\z@ \dg@ZTEMP=\tw@
   \else \divide\dg@ZTEMP\dg@YGRID \fi
   \ifnum\dg@ZTEMP>\f@ur \dg@ZTEMP=\f@ur \fi
   \ifnum\dg@ZTEMP<\@ne \dg@ZTEMP=\@ne \fi
   \unitlength=2sp\relax
   \ifnum\dg@ZTEMP<\tw@
      \advance\dg@ZTEMP\@ne
      \multiply\unitlength\dg@YGRID
   \else
      \multiply\unitlength\dg@XGRID \divide\unitlength\dg@ZTEMP
   \fi
   \dg@XGRID=\dg@ZTEMP \dg@YGRID=\tw@
   \dg@rmcommondiv\tw@\dg@XGRID\dg@YGRID
   \divide\unitlength\dg@YGRID \divide\unitlength\@m\relax}
\@namedef{dgo@..}{\let\dg@VECTOR=\dg@dotvector}%

\def\dg@dotvector(#1,#2)#3{%
   \begingroup
   \dg@XTEMP=#1\relax \dg@YTEMP=#2\relax
   \let\dg@NDOTS=\dg@XEND \let\dg@DOTDIAM=\dg@WEND
   \dg@NDOTS=\unitlength \multiply\dg@NDOTS #3\relax
   \dg@ZTEMP=\dg@YTEMP \dg@changesign\dg@YTEMP\dg@ZTEMP
   \ifnum\dg@XTEMP>\z@
      \ifnum\dg@YTEMP>\dg@XTEMP
         \multiply\dg@NDOTS\dg@YTEMP \divide\dg@NDOTS\dg@XTEMP \fi
   \else\ifnum\dg@XTEMP<\z@
      \ifnum\dg@YTEMP>-\dg@XTEMP
         \multiply\dg@NDOTS\dg@YTEMP \divide\dg@NDOTS-\dg@XTEMP \fi
   \fi\fi
   \dg@YTEMP=\dg@ZTEMP
   \divide\dg@NDOTS\dgDOTSPACING
   \ifnum\dg@NDOTS>\z@\else \dg@NDOTS=\@ne \fi
   \dg@ZTEMP=\unitlength \multiply\dg@ZTEMP #3\relax
   \divide\dg@ZTEMP\dg@NDOTS
   \ifnum\dg@XTEMP=\z@
      \dg@changesign\dg@ZTEMP\dg@YTEMP \dg@YTEMP=\dg@ZTEMP
   \else
      \dg@changesign\dg@ZTEMP\dg@XTEMP
      \multiply\dg@YTEMP\dg@ZTEMP \divide\dg@YTEMP\dg@XTEMP
      \dg@XTEMP=\dg@ZTEMP
   \fi
   \divide\dg@XTEMP\unitlength \divide\dg@YTEMP\unitlength
   \begin{picture}(0,0)
      \dg@DOTDIAM=\dgDOTSIZE \divide\dg@DOTDIAM\unitlength
      \multiput(0,0)(\dg@XTEMP,\dg@YTEMP){\dg@NDOTS}{%
         \circle*{\dg@DOTDIAM}}%
      \multiply\dg@XTEMP\dg@NDOTS \multiply\dg@YTEMP\dg@NDOTS
      \put(\dg@XTEMP,\dg@YTEMP){\vector(#1,#2){0}}%
   \end{picture}%
   \endgroup}%
\newif\ifdg@LATEXGEOM \dg@LATEXGEOMfalse

\@ifundefined{lamsvector}{}{%
   \@namedef{dgo@}{%
      \let\dg@VECTOR=\lamsvector
      \lamsshaft{?}\lamssource{?}\lamstarget{?}%
      \dg@LBLPOS=\dgARROWPARTS \divide\dg@LBLPOS\tw@\relax}%
   \@namedef{dgo@..}{\lamsshaft{.}}%
   \@namedef{dgo@=}{\lamsshaft{=}}%
   \@namedef{dgo@-}{\lamstarget{0}}%
   \@namedef{dgo@'}{\lamstarget{'}}%
   \@namedef{dgo@`}{\lamstarget{`}}%
   \@namedef{dgo@A}{\lamstarget{A}}%
   \@namedef{dgo@V}{\lamssource{V}}%
   \@namedef{dgo@J}{\lamssource{J}}%
   \@namedef{dgo@L}{\lamssource{L}}%
   \@namedef{dgo@S}{\lamssource{S}}%
   \def\dggeometry{
      \dg@ZTEMP=\dg@XGRID \multiply\dg@ZTEMP\tw@
      \ifnum\dg@YGRID=\z@ \dg@ZTEMP=\tw@
      \else \divide\dg@ZTEMP\dg@YGRID \fi
      \ifnum\dg@ZTEMP>6\relax \dg@ZTEMP=6\relax \fi
      \ifdg@LATEXGEOM\ifnum\dg@ZTEMP>\f@ur \dg@ZTEMP=\f@ur \fi\fi
      \ifnum\dg@ZTEMP<\@ne \dg@ZTEMP=\@ne \fi
      \unitlength=2sp\relax
      \ifnum\dg@ZTEMP<\tw@
         \advance\dg@ZTEMP\@ne
         \multiply\unitlength\dg@YGRID
      \else
         \multiply\unitlength\dg@XGRID \divide\unitlength\dg@ZTEMP
      \fi
      \dg@XGRID=\dg@ZTEMP \dg@YGRID=\tw@
      \dg@rmcommondiv\tw@\dg@XGRID\dg@YGRID
      \divide\unitlength\dg@YGRID \divide\unitlength\@m
      \dg@LATEXGEOMfalse}
   \def\dgt@nee{\dg@DX=\tw@ \dg@DY=\@ne \dg@SIZE=\tw@}%
   \def\dgt@see{\dg@DX=\tw@ \dg@DY=\m@ne \dg@SIZE=\tw@}%
   \def\dgt@nww{\dg@DX=-\tw@ \dg@DY=\@ne \dg@SIZE=\tw@}%
   \def\dgt@sww{\dg@DX=-\tw@ \dg@DY=\m@ne \dg@SIZE=\tw@}%
   \def\dgt@nnne{\dg@DX=\@ne \dg@DY=\thr@@ \dg@SIZE=\@ne}%
   \def\dgt@nnnw{\dg@DX=\m@ne \dg@DY=\thr@@ \dg@SIZE=\@ne}%
   \def\dgt@sssw{\dg@DX=\m@ne \dg@DY=-\thr@@ \dg@SIZE=\@ne}%
   \def\dgt@ssse{\dg@DX=\@ne \dg@DY=-\thr@@ \dg@SIZE=\@ne}%
   \def\dgt@nnnee{\dg@DX=\tw@ \dg@DY=\thr@@ \dg@SIZE=\tw@}%
   \def\dgt@nnnww{\dg@DX=-\tw@ \dg@DY=\thr@@ \dg@SIZE=\tw@}%
   \def\dgt@sssww{\dg@DX=-\tw@ \dg@DY=-\thr@@ \dg@SIZE=\tw@}%
   \def\dgt@sssee{\dg@DX=\tw@ \dg@DY=-\thr@@ \dg@SIZE=\tw@}%
   \def\dgt@nneee{\dg@DX=\thr@@ \dg@DY=\tw@ \dg@SIZE=\thr@@}%
   \def\dgt@nnwww{\dg@DX=-\thr@@ \dg@DY=\tw@ \dg@SIZE=\thr@@}%
   \def\dgt@sswww{\dg@DX=-\thr@@ \dg@DY=-\tw@ \dg@SIZE=\thr@@}%
   \def\dgt@sseee{\dg@DX=\thr@@ \dg@DY=-\tw@ \dg@SIZE=\thr@@}%
   \def\dgt@neee{\dg@DX=\thr@@ \dg@DY=\@ne \dg@SIZE=\thr@@
      \global\dg@LATEXGEOMtrue}%
   \def\dgt@nwww{\dg@DX=-\thr@@ \dg@DY=\@ne \dg@SIZE=\thr@@
      \global\dg@LATEXGEOMtrue}%
   \def\dgt@swww{\dg@DX=-\thr@@ \dg@DY=\m@ne \dg@SIZE=\thr@@
      \global\dg@LATEXGEOMtrue}%
   \def\dgt@seee{\dg@DX=\thr@@ \dg@DY=\m@ne \dg@SIZE=\thr@@
      \global\dg@LATEXGEOMtrue}%
}
\newcount\dg@HORIZ      \newcount\dg@VERT
\newcount\dg@XLPAD      \newcount\dg@YBPAD
\newcount\dg@XRPAD      \newcount\dg@YTPAD
\newcount\dg@X          \newcount\dg@Y
\newcount\dg@XGRID      \newcount\dg@YGRID
\newcount\dg@SIZE
\newcount\dg@USERSIZE
\newcount\dg@DX         \newcount\dg@DY
\newcount\dg@XLBL       \newcount\dg@YLBL
\newcount\dg@XOFFSET    \newcount\dg@YOFFSET
\newcount\dg@LBLOFF
\newcount\dg@XLBLOFF    \newcount\dg@YLBLOFF
\newcount\dg@LBLPOS
\newcount\dg@XTEMP      \newcount\dg@YTEMP
\newcount\dg@ZTEMP
\newcount\dg@XNODE      \newcount\dg@YNODE
\newcount\dg@XEND       \newcount\dg@YEND
\newcount\dg@WEND       \newcount\dg@HEND
\newcount\dg@COUNT
\newbox\dg@NODEBOX
\newoptenvironment*{diagram}{}[1]{%
   \global\advance\dg@COUNT\@ne \typeout{[diagram \the\dg@COUNT:}%
   \let\node=\dg@node \let\\=\dg@cr \let\arrow=\dg@arrow
   \def\dg@BIGNODE{#1}%
   \ifx\dg@BIGNODE\@empty\else
      \setbox\dg@NODEBOX=\hbox{$\dgeverynode{#1}$}%
      \dg@XTEMP=\wd\dg@NODEBOX
      \dg@YTEMP=\ht\dg@NODEBOX \advance\dg@YTEMP\dp\dg@NODEBOX
      \ifnum\dg@YTEMP=\z@ \dg@ZTEMP=\z@
      \else \dg@ZTEMP=\dg@XTEMP \divide\dg@ZTEMP\dg@YTEMP\fi
      \ifnum\dg@ZTEMP>\dgMAXSQUARE
         \ifnum\dg@ZTEMP<\dgMINDBLSQ \dg@XGRID=\thr@@ \dg@YGRID=\tw@
         \else \dg@XGRID=\tw@ \dg@YGRID=\@ne \fi
      \else \dg@XGRID=\@ne \dg@YGRID=\@ne \fi
      \advance\dg@XTEMP\dgHORIZPAD \advance\dg@YTEMP\dgVERTPAD
      \dg@XLPAD=\dg@XTEMP \divide\dg@XLPAD\tw@
      \dg@XRPAD=\dg@XTEMP \divide\dg@XRPAD\tw@
      \dg@YBPAD=\dg@YTEMP \divide\dg@YBPAD\tw@
      \dg@YTPAD=\dg@YTEMP \divide\dg@YTPAD\tw@
      \advance\dg@XTEMP\dgARROWLENGTH
      \divide\dg@XTEMP\dg@XGRID \divide\dg@XTEMP\@m
      \unitlength=1sp\relax \multiply\unitlength\dg@XTEMP
   \fi
   \typeout{saving...}%
   \dg@X=-\@ne \dg@Y=\z@ \dg@HORIZ=\z@\relax
   \let\dg@SLIST=\@empty
   \let\dg@NLIST=\@empty \let\dg@ALIST=\@empty
   \let\dg@PASS=\dg@savepass
}{
   \dg@VERT=-\dg@Y\relax
   \typeout{calculating...}%
   \ifx\dg@BIGNODE\@empty
      \dg@XGRID=\z@ \dg@YGRID=\z@
      \dg@XLPAD=\z@ \dg@XRPAD=\z@ \dg@YBPAD=\z@ \dg@YTPAD=\z@
      \let\dg@PASS=\dg@geompass
      \dg@NLIST \dg@ALIST
      \ifnum\dg@XGRID>\z@\else \dg@XGRID=\quad\relax \fi
      \ifnum\dg@YGRID>\z@\else \dg@YGRID=\quad\relax \fi
      \dggeometry
      \ifdim\unitlength>\z@\else \unitlength=\quad \fi
      \ifnum\dg@XGRID>\z@\else \dg@XGRID=\@ne \fi
      \ifnum\dg@YGRID>\z@\else \dg@YGRID=\@ne \fi
   \fi
   \dg@LBLOFF=\dgLABELOFFSET \divide\dg@LBLOFF\unitlength
   \multiply\dg@HORIZ\@m \multiply\dg@HORIZ\dg@XGRID
   \multiply\dg@VERT\@m \multiply\dg@VERT\dg@YGRID
   \divide\dg@XLPAD\unitlength \divide\dg@XRPAD\unitlength
   \divide\dg@YBPAD\unitlength \divide\dg@YTPAD\unitlength
   \typeout{drawing...}%
   \let\dg@PASS=\dg@drawpass
   \hspace*{\dg@XLPAD\unitlength}%
   \vcenter{%
      \hsize=0pt\relax\parindent=0pt\relax
      \parskip=0pt\relax\baselineskip=0pt\relax
      \vspace*{\dg@YTPAD\unitlength}%
      \begin{picture}(0,0)(0,0)\dg@NLIST\dg@ALIST\end{picture}%
      \raisebox{\z@}[\z@][\dg@VERT\unitlength]{}%
      \vspace*{\dg@YBPAD\unitlength}%
   }
   \hspace*{\dg@HORIZ\unitlength}%
   \hspace*{\dg@XRPAD\unitlength}%
   \typeout{done]}}%

\def\dg@savepass{s}
\def\dg@geompass{g}
\def\dg@drawpass{d}

\newoptcommand{\dg@node}{\@ne}[2]{%
   \ifx\dg@PASS\dg@savepass
      \dg@XTEMP=#1\relax
      \ifnum\dg@XTEMP<\@ne \dg@XTEMP=\@ne\fi
      \advance\dg@X\dg@XTEMP
      \ifnum\dg@HORIZ<\dg@X \dg@HORIZ=\dg@X \fi
      \setbox\dg@NODEBOX=\hbox{$\dgeverynode{#2}$}%
      \dg@XTEMP=\wd\dg@NODEBOX \advance\dg@XTEMP\dgHORIZPAD
      \dg@YTEMP=\ht\dg@NODEBOX \advance\dg@YTEMP\dp\dg@NODEBOX
      \advance\dg@YTEMP\dgVERTPAD
      \toks\z@=\expandafter{\dg@SLIST}%
      \edef\dg@SLIST{\the\toks\z@
         ,\noexpand\dg@XNODE=\number\dg@X\noexpand\relax
         \noexpand\dg@YNODE=\number\dg@Y\noexpand\relax
         \noexpand\dg@XTEMP=\number\dg@XTEMP\noexpand\relax
         \noexpand\dg@YTEMP=\number\dg@YTEMP\noexpand\relax}%
      \toks\z@=\expandafter{\dg@NLIST}%
      \toks\tw@={\dg@node{#2}}%
      \edef\dg@NLIST{\the\toks\z@
         \noexpand\dg@X=\number\dg@X\noexpand\relax
         \noexpand\dg@Y=\number\dg@Y\noexpand\relax
         \the\toks\tw@}%
   \else\ifx\dg@PASS\dg@geompass
      \ifnum\dg@X=\z@
         \dg@getnodesize
            {\dg@SLIST}{\dg@X}{\dg@Y}{\dg@WEND}{\dg@HEND}%
         \divide\dg@WEND\tw@
         \ifnum\dg@XLPAD<\dg@WEND \dg@XLPAD=\dg@WEND \fi\fi
      \ifnum\dg@X=\dg@HORIZ
         \dg@getnodesize
            {\dg@SLIST}{\dg@X}{\dg@Y}{\dg@WEND}{\dg@HEND}%
         \divide\dg@WEND\tw@
         \ifnum\dg@XRPAD<\dg@WEND \dg@XRPAD=\dg@WEND \fi\fi
      \ifnum\dg@Y=\z@
         \dg@getnodesize
            {\dg@SLIST}{\dg@X}{\dg@Y}{\dg@WEND}{\dg@HEND}%
         \divide\dg@HEND\tw@
         \ifnum\dg@YTPAD<\dg@HEND \dg@YTPAD=\dg@HEND \fi\fi
      \ifnum\dg@Y=-\dg@VERT\relax
         \dg@getnodesize
            {\dg@SLIST}{\dg@X}{\dg@Y}{\dg@WEND}{\dg@HEND}%
         \divide\dg@HEND\tw@
         \ifnum\dg@YBPAD<\dg@HEND \dg@YBPAD=\dg@HEND \fi\fi
   \else\ifx\dg@PASS\dg@drawpass
      \dg@XNODE=\dg@X \multiply\dg@XNODE\@m
      \multiply\dg@XNODE\dg@XGRID
      \dg@YNODE=\dg@Y \multiply\dg@YNODE\@m
      \multiply\dg@YNODE\dg@YGRID
      \setbox\dg@NODEBOX=\hbox{$\dgeverynode{#2}$}%
      \put(\dg@XNODE,\dg@YNODE){\dg@makebox{\box\dg@NODEBOX}}%
   \fi\fi\fi}%

\newoptcommand{\dg@cr}{\@ne}[1]{%
   \ifx\dg@PASS\dg@savepass
      \dg@YTEMP=#1\relax
      \ifnum\dg@YTEMP<\@ne \dg@YTEMP=\@ne \fi
      \advance\dg@Y -\dg@YTEMP\relax
      \dg@X=-\@ne\relax\fi}%
\newoptcommand{\dg@arrow}{\@ne}[2]{%
   \begingroup
   \dg@USERSIZE=#1\relax
   \ifnum\dg@USERSIZE<\@ne \dg@USERSIZE=\@ne \fi
   \dg@parse{#2}%
   \ifx\dg@PASS\dg@savepass
      \ifx\dg@label\dgl@b \let\dg@label=\dgl@t \fi
      \ifx\dg@label\dgl@r \let\dg@label=\dgl@l \fi
      \let\dg@process=\dg@save
   \else\ifx\dg@PASS\dg@geompass
      \let\dg@process=\dg@ignore
      \dg@geomcalc
   \else\ifx\dg@PASS\dg@drawpass
      \let\dg@process=\dg@draw
      \dg@drawcalc
   \fi\fi\fi
   \dg@label{\dg@process{#1}{#2}}}%

\newoptcommand{\arrow}{\@ne}[2]{%
   \dg@parse{#2}%
   \ifx\dg@label\dgl@b \let\dg@label=\dgl@t \fi
   \ifx\dg@label\dgl@r \let\dg@label=\dgl@l \fi
   \dg@label{\dg@textarrow{#1}{#2}}}%

\def\dg@textarrow#1#2#3#4{%
   \mathop{{\dgHORIZPAD=0pt\relax\dgVERTPAD=0pt\relax
      \begin{diagram}
         \node{}\arrow[#1]{#2}{#3}{#4}\node{}
      \end{diagram}}}}

\def\dg@parse#1{%
   \let\dg@label=\dgl@ \dgo@
   \let\dg@type=\@empty \let\dg@lbltype=\@empty
   \@for\dg@list:=#1\do{%
      \ifx\dg@type\@empty \let\dg@type=\dg@list
      \else\ifx\dg@lbltype\@empty \let\dg@lbltype=\dg@list
         \@ifundefined{dgo@\dg@list}{}{\@nameuse{dgo@\dg@list}}%
      \else
         \@ifundefined{dgo@\dg@list}{}{\@nameuse{dgo@\dg@list}}%
      \fi\fi}%
   \@ifundefined{dgt@\dg@type}{\dgt@e}{\@nameuse{dgt@\dg@type}}%
   \@ifundefined{dgl@\dg@lbltype}{}{%
      \dg@letname\dg@label{dgl@\dg@lbltype}}}

\def\dg@draw#1#2#3#4{%
   \put(\dg@X,\dg@Y){\dg@makebox{%
      \begin{picture}(0,0)%
         \thinlines
         \put(\dg@XOFFSET,\dg@YOFFSET){%
            \dg@VECTOR(\dg@DX,\dg@DY){\dg@SIZE}}%
         \put(\dg@XLBL,\dg@YLBL){\dg@makebox{%
            \begin{picture}(0,0)%
               \put(\dg@XLBLOFF,\dg@YLBLOFF){%
                  \dg@makebox[\dg@LBLONE]{$\dgeverylabel{#3}$}}%
               \put(-\dg@XLBLOFF,-\dg@YLBLOFF){%
                  \dg@makebox[\dg@LBLTWO]{$\dgeverylabel{#4}$}}%
            \end{picture}}}%
      \end{picture}}}%
   \endgroup}%

\def\dg@save#1#2#3#4{%
   \endgroup 
   \toks\z@=\expandafter{\dg@ALIST}%
   \toks\tw@={\dg@arrow[#1]{#2}{#3}{#4}}%
   \edef\dg@ALIST{\the\toks\z@%
      \noexpand\dg@X=\number\dg@X\noexpand\relax
      \noexpand\dg@Y=\number\dg@Y\noexpand\relax
      \the\toks\tw@}}%

\def\dg@ignore#1#2#3#4{\endgroup}

\def\dg@geomcalc{%
   \dg@XEND=\dg@SIZE \multiply\dg@XEND\dg@USERSIZE
   \ifnum\dg@DX=\z@
      \dg@YEND=\dg@XEND \dg@XEND=\z@
      \dg@changesign\dg@YEND\dg@DY
   \else
      \dg@changesign\dg@XEND\dg@DX \dg@YEND=\dg@XEND
      \multiply\dg@YEND\dg@DY \divide\dg@YEND\dg@DX
   \fi
   \advance\dg@XEND\dg@X \advance\dg@YEND\dg@Y
   \dg@getnodesize
      {\dg@SLIST}{\dg@XEND}{\dg@YEND}{\dg@WEND}{\dg@HEND}%
   \dg@XOFFSET=\dg@WEND \dg@YOFFSET=\dg@HEND
   \dg@getnodesize
      {\dg@SLIST}{\dg@X}{\dg@Y}{\dg@WEND}{\dg@HEND}%
   \advance\dg@XOFFSET\dg@WEND \divide\dg@XOFFSET\tw@
   \advance\dg@YOFFSET\dg@HEND \divide\dg@YOFFSET\tw@
   \dg@XTEMP=\dgARROWLENGTH \dg@YTEMP=\dgARROWLENGTH
   \ifnum\dg@DX>\z@
      \dg@ZTEMP=\dg@DX \multiply\dg@XTEMP\dg@DX
   \else \dg@ZTEMP=-\dg@DX \multiply\dg@XTEMP -\dg@DX \fi
   \ifnum\dg@DY>\z@
      \advance\dg@ZTEMP\dg@DY \multiply\dg@YTEMP\dg@DY
   \else \advance\dg@ZTEMP -\dg@DY \multiply\dg@YTEMP -\dg@DY\fi
   \ifnum\dg@ZTEMP=\z@\else
      \divide\dg@XTEMP\dg@ZTEMP \divide\dg@YTEMP\dg@ZTEMP
      \advance\dg@XOFFSET\dg@XTEMP \advance\dg@YOFFSET\dg@YTEMP
   \fi
   \divide\dg@XOFFSET\dg@SIZE \divide\dg@XOFFSET\dg@USERSIZE
   \divide\dg@YOFFSET\dg@SIZE \divide\dg@YOFFSET\dg@USERSIZE
   \ifnum\dg@DX=\z@ \dg@XOFFSET=\z@ \fi
   \ifnum\dg@DY=\z@ \dg@YOFFSET=\z@ \fi
   \ifnum\dg@XGRID<\dg@XOFFSET \global\dg@XGRID=\dg@XOFFSET\fi
   \ifnum\dg@YGRID<\dg@YOFFSET \global\dg@YGRID=\dg@YOFFSET\fi
   \relax}

\def\dg@drawcalc{%
   \dg@XEND=\dg@SIZE \multiply\dg@XEND\dg@USERSIZE
   \ifnum\dg@DX=\z@
      \dg@YEND=\dg@XEND \dg@XEND=\z@
      \dg@changesign\dg@YEND\dg@DY
   \else
      \dg@changesign\dg@XEND\dg@DX \dg@YEND=\dg@XEND
      \multiply\dg@YEND\dg@DY \divide\dg@YEND\dg@DX
   \fi
   \advance\dg@XEND\dg@X \advance\dg@YEND\dg@Y
   \dg@getnodesize
      {\dg@SLIST}{\dg@XEND}{\dg@YEND}{\dg@WEND}{\dg@HEND}%
   \divide\dg@WEND\unitlength \divide\dg@HEND\unitlength
   \multiply\dg@DX\dg@XGRID \multiply\dg@DY\dg@YGRID
   \dg@rmcommondiv\tw@\dg@DX\dg@DY
   \dg@rmcommondiv\tw@\dg@DX\dg@DY 
   \dg@rmcommondiv\thr@@\dg@DX\dg@DY
   \multiply\dg@SIZE\dg@USERSIZE \multiply\dg@SIZE\@m
   \ifnum\dg@DX=\z@
      \multiply\dg@SIZE\dg@YGRID
      \divide\dg@HEND\tw@ \advance\dg@SIZE -\dg@HEND
      \dg@getnodesize
         {\dg@SLIST}{\dg@X}{\dg@Y}{\dg@WEND}{\dg@YOFFSET}%
      \divide\dg@YOFFSET\unitlength \divide\dg@YOFFSET\tw@
      \advance\dg@SIZE -\dg@YOFFSET
      \dg@XOFFSET=\z@
      \def\dg@LBLONE{r}\def\dg@LBLTWO{l}%
      \dg@XLBL=\z@ \dg@YLBL=\dg@SIZE
      \multiply\dg@YLBL\dg@LBLPOS
      \divide\dg@YLBL\dgARROWPARTS\relax
      \advance\dg@YLBL\dg@YOFFSET
      \dg@changesign\dg@YLBL\dg@DY
      \dg@changesign\dg@YOFFSET\dg@DY
   \else
      \multiply\dg@SIZE\dg@XGRID
      \ifnum\dg@DY=\z@
         \divide\dg@WEND\tw@ \advance\dg@SIZE -\dg@WEND
         \dg@getnodesize
            {\dg@SLIST}{\dg@X}{\dg@Y}{\dg@XOFFSET}{\dg@HEND}%
         \divide\dg@XOFFSET\unitlength \divide\dg@XOFFSET\tw@
         \advance\dg@SIZE -\dg@XOFFSET
         \dg@YOFFSET=\z@
         \def\dg@LBLONE{b}\def\dg@LBLTWO{t}%
         \dg@YLBL=\z@ \dg@XLBL=\dg@SIZE
         \multiply\dg@XLBL\dg@LBLPOS
         \divide\dg@XLBL\dgARROWPARTS\relax
         \advance\dg@XLBL\dg@XOFFSET
         \dg@changesign\dg@XLBL\dg@DX
         \dg@changesign\dg@XOFFSET\dg@DX
      \else
         \divide\dg@WEND\tw@ \divide\dg@HEND\tw@
         \multiply\dg@HEND\dg@DX \divide\dg@HEND\dg@DY
         \ifnum\dg@HEND<\z@ \multiply\dg@HEND\m@ne \fi
         \ifnum\dg@WEND<\dg@HEND \advance\dg@SIZE -\dg@WEND
         \else \advance\dg@SIZE -\dg@HEND \fi
         \dg@getnodesize
            {\dg@SLIST}{\dg@X}{\dg@Y}{\dg@WEND}{\dg@HEND}%
         \divide\dg@WEND\unitlength \divide\dg@WEND\tw@
         \divide\dg@HEND\unitlength \divide\dg@HEND\tw@
         \multiply\dg@HEND\dg@DX \divide\dg@HEND\dg@DY
         \ifnum\dg@HEND<\z@ \multiply\dg@HEND\m@ne \fi
         \ifnum\dg@WEND<\dg@HEND \dg@XOFFSET=\dg@WEND
         \else \dg@XOFFSET=\dg@HEND \fi
         \advance\dg@SIZE -\dg@XOFFSET
         \dg@changesign\dg@XOFFSET\dg@DX
         \dg@YOFFSET=\dg@XOFFSET
         \multiply\dg@YOFFSET\dg@DY \divide\dg@YOFFSET\dg@DX
         \def\dg@LBLONE{br}\def\dg@LBLTWO{tl}%
         \ifnum\dg@DX<\z@ \ifnum\dg@DY>\z@
            \def\dg@LBLONE{bl}\def\dg@LBLTWO{tr}\fi\fi
         \ifnum\dg@DX>\z@ \ifnum\dg@DY<\z@
            \def\dg@LBLONE{bl}\def\dg@LBLTWO{tr}\fi\fi
         \dg@XLBL=\dg@SIZE
         \multiply\dg@XLBL\dg@LBLPOS
         \divide\dg@XLBL\dgARROWPARTS\relax
         \dg@changesign\dg@XLBL\dg@DX
         \dg@YLBL=\dg@XLBL
         \multiply\dg@YLBL\dg@DY \divide\dg@YLBL\dg@DX
         \advance\dg@XLBL\dg@XOFFSET
         \advance\dg@YLBL\dg@YOFFSET
      \fi
   \fi
   \dg@XLBLOFF=-\dg@DY \dg@changesign\dg@XLBLOFF\dg@DX
   \dg@YLBLOFF=\dg@DX \dg@changesign\dg@YLBLOFF\dg@DX
   \ifnum\dg@DX=\z@ \dg@XLBLOFF=\m@ne \fi
   \dg@XTEMP=\dg@DX \dg@changesign\dg@XTEMP\dg@DX
   \dg@YTEMP=\dg@DY \dg@changesign\dg@YTEMP\dg@DY
   \ifnum\dg@YTEMP>\dg@XTEMP \dg@XTEMP=\dg@YTEMP \fi
   \ifnum\dg@XTEMP=\z@ \dg@XTEMP=\@ne \fi
   \multiply\dg@XLBLOFF\dg@LBLOFF \divide\dg@XLBLOFF\dg@XTEMP
   \multiply\dg@YLBLOFF\dg@LBLOFF \divide\dg@YLBLOFF\dg@XTEMP
   \multiply\dg@X\@m \multiply\dg@X\dg@XGRID
   \multiply\dg@Y\@m \multiply\dg@Y\dg@YGRID
   \relax}%

\def\dg@rmcommondiv#1#2#3{%
   \dg@XTEMP=#2\relax
   \divide\dg@XTEMP #1\relax \multiply\dg@XTEMP #1\relax
   \dg@YTEMP=#3\relax
   \divide\dg@YTEMP #1\relax \multiply\dg@YTEMP #1\relax
   \ifnum\dg@XTEMP=#2\relax \ifnum\dg@YTEMP=#3\relax
      \divide#2#1\relax \divide#3#1\relax \fi\fi}%

\def\dg@changesign#1#2{%
   \ifnum #2<\z@ \multiply#1\m@ne
   \else\ifnum #2=\z@ #1=\z@ \fi\fi}%

\def\dg@getnodesize#1#2#3#4#5{%
   #4=\z@\relax #5=\z@\relax
   \expandafter\@for\expandafter\dg@trynode
   \expandafter:\expandafter=#1\do{%
      \dg@XNODE=\m@ne 
      \dg@trynode
      \ifnum #2=\dg@XNODE \ifnum #3=\dg@YNODE
         #4=\dg@XTEMP\relax #5=\dg@YTEMP\relax\fi\fi}}%

\newoptcommand{\dg@makebox}{}[2]{%
   \expandafter\makebox\expandafter(\expandafter
      0\expandafter,\expandafter0\expandafter)\expandafter
      [#1]{#2}}%

\def\dg@novector(#1,#2)#3{}%

\def\dg@letname#1#2{%
   \relax\expandafter
   \let\expandafter #1\csname #2\endcsname\relax}%

\def\dgl@#1{#1{}{}}%
\def\dgl@t#1#2{#1{#2}{}}%
\def\dgl@b#1#2{#1{}{#2}}%
\def\dgl@tb#1#2#3{#1{#2}{#3}}%
\def\dgl@l#1#2{#1{#2}{}}%
\def\dgl@r#1#2{#1{}{#2}}%
\def\dgl@lr#1#2#3{#1{#2}{#3}}%
\makeatother
\newcommand{\jac}{{\frak J}}
\newcommand{\mod}{{\frak M}^0_g}

\newcommand{\univ}{{\frak C}_g}
\newcommand{\symm}{{\frak C}_g^{(d)}}

\newcommand{\lra}{\longrightarrow}

\newcommand{\pr}{{\Bbb P}}

\newcommand{\hilb}{{\frak H}_{d,g,r} }
\newcommand{\unhilb}{{\frak C}_{\frak H}}

\newtheorem{lem}{Lemma}[section]

\newtheorem{theor}{Theorem}[section]
\newtheorem{defi}{Definition}[section]
\newtheorem{rem}{Remark}[section]

\begin{document}
%
\expandafter\ifx\csname amssym.def\endcsname\relax \else\endinput\fi
%
\expandafter\edef\csname amssym.def\endcsname{%
       \catcode`\noexpand\@=\the\catcode`\@\space}
\catcode`\@=11
%

\def\undefine#1{\let#1\undefined}
\def\newsymbol#1#2#3#4#5{\let\next@\relax
 \ifnum#2=\@ne\let\next@\msafam@\else
 \ifnum#2=\tw@\let\next@\msbfam@\fi\fi
 \mathchardef#1="#3\next@#4#5}
\def\mathhexbox@#1#2#3{\relax
 \ifmmode\mathpalette{}{\m@th\mathchar"#1#2#3}%
 \else\leavevmode\hbox{$\m@th\mathchar"#1#2#3$}\fi}
\def\hexnumber@#1{\ifcase#1 0\or 1\or 2\or 3\or 4\or 5\or 6\or 7\or
8\or
 9\or A\or B\or C\or D\or E\or F\fi}

\font\tenmsa=msam10
\font\sevenmsa=msam7
\font\fivemsa=msam5
\newfam\msafam
\textfont\msafam=\tenmsa
\scriptfont\msafam=\sevenmsa
\scriptscriptfont\msafam=\fivemsa
\edef\msafam@{\hexnumber@\msafam}
\mathchardef\dabar@"0\msafam@39
\def\dashrightarrow{\mathrel{\dabar@\dabar@\mathchar"0\msafam@4B}}
\def\dashleftarrow{\mathrel{\mathchar"0\msafam@4C\dabar@\dabar@}}
\let\dasharrow\dashrightarrow
\def\ulcorner{\delimiter"4\msafam@70\msafam@70 }
\def\urcorner{\delimiter"5\msafam@71\msafam@71 }
\def\llcorner{\delimiter"4\msafam@78\msafam@78 }
\def\lrcorner{\delimiter"5\msafam@79\msafam@79 }
\def\yen{{\mathhexbox@\msafam@55 }}
\def\checkmark{{\mathhexbox@\msafam@58 }}
\def\circledR{{\mathhexbox@\msafam@72 }}
\def\maltese{{\mathhexbox@\msafam@7A }}

\font\tenmsb=msbm10
\font\sevenmsb=msbm7
\font\fivemsb=msbm5
\newfam\msbfam
\textfont\msbfam=\tenmsb
\scriptfont\msbfam=\sevenmsb
\scriptscriptfont\msbfam=\fivemsb
\edef\msbfam@{\hexnumber@\msbfam}
\def\Bbb#1{{\fam\msbfam\relax#1}}
\def\widehat#1{\setbox\z@\hbox{$\m@th#1$}%
 \ifdim\wd\z@>\tw@ em\mathaccent"0\msbfam@5B{#1}%
 \else\mathaccent"0362{#1}\fi}
\def\widetilde#1{\setbox\z@\hbox{$\m@th#1$}%
 \ifdim\wd\z@>\tw@ em\mathaccent"0\msbfam@5D{#1}%
 \else\mathaccent"0365{#1}\fi}
\font\teneufm=eufm10
\font\seveneufm=eufm7
\font\fiveeufm=eufm5
\newfam\eufmfam
\textfont\eufmfam=\teneufm
\scriptfont\eufmfam=\seveneufm
\scriptscriptfont\eufmfam=\fiveeufm
\def\frak#1{{\fam\eufmfam\relax#1}}
\let\goth\frak

\csname amssym.def\endcsname

\addtolength{\baselineskip}{5pt}
$ $
 \vskip.4in
 \noindent
\begin{center}
{\bf {Picard groups of Hilbert schemes of Curves}} \\
  \vskip.1in
 Alexis Kouvidakis\footnote{New address: Department of Mathematics,
 University of Crete,  Iraklion 71409, Greece.} \\
 University of Pennsylvania \\
 Department of Mathematics, DRL \\
 Philadelphia, PA 19104-6395\\
 e-mail: alex@math.upenn.edu\\
 \end{center}
 \vskip.3in
\noindent
{\bf {Abstract:}} {\em {We calculate the Picard group, over the
integers, of
the Hilbert scheme of smooth, irreducible, non-degenerate curves of
degree $d$
and  genus
 $g \geq 4$  in ${\Bbb P}^r$, in the case when $d \geq 2g+1 $ and $r
\leq d-g$.
  We express the classes of the generators in terms of some
``natural'' divisor
classes.}}
\begin{center}
 {\bf Notation and conventions}
\end{center}
\noindent
$\mod $ : Moduli space of smooth, irreducible  curves of genus $g
\geq 4$,
without automorphisms. \\
$\pi :\, \univ  \lra \mod $ : Universal curve over $\mod $. \\
$J^d(C)$ : Jacobian variety which parametrizes line bundles of degree
$d$ on
the curve $C$.\\
$\psi : \, \jac ^d\lra \mod $ : Universal Jacobian variety;  the
fiber over
$[C] \in \mod $ is $J^d(C)$. \\
$  \hilb $ : Hilbert scheme of  smooth, irreducible, non-degenerate
curves of
degree $d$
 and  genus $g$  in $\pr ^r$. \\
$q : \, \unhilb \lra \hilb $ : Universal curve over $\hilb $.\\
${\cal F} $ : Tautological  line bundle over $\unhilb $. \\
$\pr (V)$ : denotes the space of  one dimensional subspaces of $V$.\\

\section{Introduction}
Given positive integers $d,g,r$, let $\hilb $ denote the Hilbert
scheme of
smooth, irreducible, non-degenerate  curves
of degree $d$ and genus $g$ in ${\Bbb P}^r$. In general,  the
geometric
structure
 of $\hilb $ is difficult to describe and in
most of the cases is unknown but for $d \geq  2g+1$ and $r \leq
d-g$, it turns
out that $\hilb $ is smooth and irreducible, see \cite{Ha},
\cite{CS}. The
natural forgetful  map  $\hilb  \lra {\frak M}_g $, where  ${\frak
M}_g $
is the moduli space of smooth curves of genus $g \geq 4$,  is onto.
The purpose
of this paper is to describe the Picard group of $\hilb $ over the
integers
when $d \geq 2g+1$ and $r \leq  d-g$. From now on, we are going to
exclude from
the Hilbert scheme those points which represent curves with
automorphisms. By
doing so, the description of the Picard group is not effected, since
the locus
of such points is of big codimension when $g \geq 4$.

Let $\mod $ denote the moduli space of smooth, irreducible curves of
genus $g
\geq 4$ without automorphisms.  Over $\mod $ we have the universal
family $\pi
: \univ \lra \mod$. To that we can associate the family $\psi : \jac
^d \lra
\mod $, the universal Picard variety of degree $d$, whose fiber over
$[C] \in
\mod $ is  $J^d(C)$.  Over $\jac ^d$ one can construct a projective
fibration
$ \phi : {\frak P}_d \lra \jac ^d   $ whose fiber over a point $[L]
\in J^d
(C)$ is $\pr ({\Bbb C}^{r+1} \otimes  H^0(C,L))$. Note that for
$r=0$, this is
just the universal symmetric product $\symm \lra \jac ^d $ of degree
$d$. The
existence of such a bundle is based on the existence of the bundle
$\symm \lra
\jac ^d $ and the  existence of  a local section (in the analytic
topology) of
the map
 $\pi  : \, \univ \lra \mod $.

 The variety $\hilb $ can be included in ${\frak P}_d$ as follows. An
element
$h$ in $\hilb$ corresponds to a smooth, irreducible, non-degenerate
curve $C$
of degree $d$ and genus $g$ in $\pr ^r$. Let $H_i, \; i=1, \ldots ,
r+1$,
denote the hyperplane section $X_i=0$. Then, define $L\stackrel{\rm
def}{=}{\cal O}(1)|_C \in J^d(C)$ and $s_i \stackrel{\rm def}{=}
H_i|_C \in
H^0(C,L)$. We then correspond to $h$ the point $[C,L, <s_1, \ldots
,s_{r+1}>]
\in {\frak P}_d$ and the map is obviously one to one.  By doing so,
one can
factor the canonical map  $\hilb  \lra \mod   $ as
 $$ \hilb \stackrel{i}{ \hookrightarrow} {\frak P}_d \stackrel{\phi
}{\lra }
\jac ^d \stackrel{\psi }{\lra } \mod  . $$
The complement of $\hilb $ in ${\frak P}_d$ corresponds to those
tuples
  $[C,L, <s_1, \ldots ,s_{r+1}>] $ for which the space of  sections
$<s_1,
\ldots ,s_{r+1}>$
 has either a base point or does not separate points and tangent
directions on
the curve or the dimension of the span$<s_1, \ldots ,s_{r+1}>  $ is
$\leq r$.
In the first case the map of $C$ to $\pr ^r$  defined by the above
data is of
degree $<d$, in the second the map is not an embedding  and in the
third the
map is degenerate.

Over $\hilb $ we have the universal curve  $q: \, \unhilb \lra \hilb
$. By
construction $\unhilb \subset \hilb \times \pr ^r$. On $\unhilb $, we
have the
tautological bundle ${\cal F}$ which is the pull back by the
projection  map
$\eta  : \unhilb \lra {\Bbb P}^r$
of ${\cal O}(1)$ on $\pr ^r$ .  If $\pi : \, \univ \lra \mod $ is the
universal
curve over $\mod $ and $\psi $ the map as above, then we denote by
${\frak
C}_{g,d}$ the fiber product  ${\frak C}_{g,d} \stackrel{\rm
def}{=}\univ \times
_{\mod } \jac ^d $. Let $\nu $ denote the projection map ${\frak
C}_{g,d}\lra
\jac ^d $. The basic diagram we are going to use is the following:

\begin{equation}
\begin{diagram}[{\frak C}_{g,d}aaa]
 \node{} \node{\unhilb} \arrow{e,t}{{ \alpha }} \arrow{s,r}{ {q} }
\arrow{sw,l}{{\eta }} \node{{\frak C}_{g,d}} \arrow{e,t}{{\beta }}
\arrow{s,r}{{\nu }} \node{\univ } \arrow{s,r}{{\pi } }\\
\node{\pr ^r}\node{\hilb} \arrow{e,t}{{\phi } } \node{\jac ^d}
\arrow{e,t}{{\psi }} \node{\mod }
\end{diagram}
\label{diag1}
\end{equation}

\vskip.2in
 Let  $\omega _q $ denote the  relative dualizing sheaf of the map
$q: \,
\unhilb \lra \hilb $. Then we have  on $\hilb $ the following three
natural
divisor classes:

\begin{tabbing}
nn \= nnnnn   \kill
\> $A = q_{*}({\cal F}^2)$, \\
\> $B  = q_{*} ({\cal F} \omega _q)$,\\
\> $C =q_{*} (\omega _q ^2)$.
 \end{tabbing}
In this paper we describe  the Picard group of $\hilb $ and we
express the
classes of its generators in terms of the classes $A, B, C$ given
above. This
is the content of Theorem \ref{main}.

\section{Intersection calculations.}
\noindent
We are going to use the following results
\begin{enumerate}
\item The Picard group of $\mod $ is freely generated over the
integers by the
determinant $\lambda $ of the Hodge bundle. In other words, $\lambda
=det \pi
_{*} \omega_{\pi }$, where $\pi :\univ  \lra \mod $ is the universal
curve, see
 \cite{AC2}.
\item   The Picard group of the universal Jacobian $\psi  :\jac
^d\lra \mod $
is freely generated over the integers by the pull back $\psi
^*\lambda$ and a
line bundle ${\cal L }_d$.  The later is  uniquely defined up to the
pull back
of line bundles  from $\mod $  and has the following  characteristic
property:
 its restriction to a fiber $J^d(C)$ has class
  $ k_d \theta $, where $k_d=\frac{2g-2}{g.c.d.(2g-2,d-g+1)}$ and
$\theta $
denotes the
   class of the theta divisor, see \cite{K1}.
\item On $\jac ^d \times _{\mod } \univ $ there is a line bundle
${\cal P }_d $
with the property ${\cal P }_d|_{[L] \times C} \cong L^{\otimes
s_d}$, where
$s_d=gcd(2g-2,d+g-1)$ is minimum with this property, see Application
in Section
5 of  \cite{K2}.
\end{enumerate}

 There are various ways to construct the bundles ${\cal L}_d$ and
${\cal P}_d$
in $2.$ and $3.$ above. In the following we give  a construction
which is the
most convenient for our purposes.

  We start with some notation.  On the $d$-th symmetric product of a
curve we
denote by $\Delta $  the diagonal or its class in the Chow ring. We
denote by
$x$ the class of the
image of a coordinate plane from the $d$-th ordinary product under
the
canonical map. Let $u: C^{(d)} \lra J^d(C)$ be  the Abel-Jacobi map
and let
$\omega _u$ the relative  dualizing sheaf.  Given a line bundle $M$
on $C$,
then we associate to that a line bundle $L_M$ on the symmetric
product
$C^{(d)}$ as follows: Take a meromorphic section of M
 written in the form $D_1-D_2$, where $D_1,\, D_2$  are effective
divisors
which are sums  of distinct points. Define on  $C^{(d)}$  the
divisors
$X_{D_i}, \; i=1,\, 2$ with support $\{D \in  C_d, \; \mbox{s.t.}\;
D \cap D_i
\neq \emptyset   \}$.  Then we define $L_M$ to
be the line bundle ${\cal O}(X_{D_1}) \otimes {\cal
O}(X_{D_2})^{-1}$. It is
easy to see that this is independent from the choice of $D_1$ and
$D_2$. For
$M=K$, we denote by $L_K$ the line bundle which corresponds to the
canonical
divisor $K$.
\begin{rem}
{\rm For $ d \geq  2g+1 $, the Abel-Jacobi map $u: C^{(d)} \lra
J^d(C)$ is a
fibration of projective spaces of dimension $r=d-g$. In that case, if
$A \in
J^d(C)$, then $L_M|_{u^{-1}(A) \simeq {\Bbb P}^{\, r}} \simeq {\cal
O}_{{\Bbb
P}^{\, r}}
 ({\rm deg}M)$.}
\label{rem21}
\end{rem}

We describe now the analogue of $L_K$   for families of curves.
Consider the
diagram

\begin{equation}
\begin{diagram}[aa]
\node{\univ ^{\times d}}  \arrow{s,l}{\pi _i}  \arrow[2]{e,t}{c}
\arrow{se,t}{\chi}  \node[2]{\univ ^{(d)}} \arrow{sw,t}{\chi ^{'}}
\\
\node{\univ}  \arrow{e,t}{\pi }   \node{\mod }
\label{diag2}
\end{diagram}
\end{equation}
\noindent
where the maps $ c, \chi  , \chi ^{'}  , \pi $ are the canonical ones
and $\pi
_i $ is the $i$-th projection.

 For $l$ big enough,  let $\sigma (l)$ be a generic section of  $
H^0(\univ ,
\omega _{\pi } \otimes  \pi ^{*}\lambda ^{\otimes l}) \cong H^0(\mod
, \pi _{*}
\omega _{\pi } \otimes \lambda ^{\otimes l}) $. Since $\pi _* \omega
_{\pi }$
is a locally free sheaf of rank $g$, we may assume that the section
$\sigma
(l)$  does not vanish identically on the fibers of $\pi $ over a
Zariski open
${\cal U}$  of  $\mod $ with complement of codimension $g$.  Then, as
above, we
can construct over ${\cal U}$ a divisor whose restriction to the
fiber
$C^{(d)}$ is the divisor $X_{\sigma (l)|_{C}}$  (following the above
notation).
 The corresponding line bundle extends uniquely to a line bundle on
$\univ
^{(d)}$.  We denote that by
 ${\cal L}_{\sigma (l)}$. On $\univ ^{\times d}$   we define  ${\cal
L}_K =
\otimes _{i=1}^d
 \pi ^{*}_i     \omega _{\pi }$. Then we have:
 \begin{lem}
Following the above notation, we have
\[
c^{*}{\cal L}_{\sigma (l)} \cong   {\cal L}_K  \otimes  \chi
^{*}(\lambda
^{\otimes dl}).
\]
 \label{pull back}
 \end{lem}
{\sc {Proof:}}  By construction we have that $c^{*}{\cal L}_{\sigma
(l)} \cong
\otimes _{i=1}^d  \pi _i ^{*} {\cal O}(\sigma(l))  \cong  \otimes
_{i=1}^d  \pi
_i ^{*} ( \omega _{\pi } \otimes \pi ^{*}\lambda ^{\otimes l})  \cong
{\cal
L}{_K} \otimes  \chi ^{*}(\lambda ^{\otimes dl}) $.
 \begin{flushright}
   $\Box$
   \end{flushright}

\begin{defi}
{\rm   We define on $\symm $ the line bundle   ${\cal L}_{\omega }
\stackrel{\rm def}{=}
  {\cal L}_{\sigma  (l)}  \otimes   \chi ^{'*}  (\lambda ^{-dl})$. }
\end{defi}

\begin{rem}
{\rm  The characteristic property of $ {\cal L}_{\omega }  $ is that
$c^{*}
 {\cal L}_{\omega }  \cong  {\cal L}_K$.  This is an application of
the above
Lemma
 \ref{pull back}.  Note also that the line bundle  ${\cal L}_{\omega
} $ does
not depend on the choice of the number $l$ and that the restriction
of ${\cal
L}_{\omega }$ to a fiber $C^{(d)}$ is isomorphic to $L_K$.}
\label{rem22}
\end{rem}

  Before we continue we need two more things: The first is  that the
universal
symmetric product carries a universal bundle ${\cal D}$. Indeed, if
$\delta :
{\frak C}^{(d-1)} \times \univ \lra \symm  \times \univ $ is the map
sending
the pair $(D,p)$ to $(D+p,p)$, then the line bundle ${\cal D}$
corresponds to
the divisor which is
the image of the above map. The second is the MacDonnald's formula
which
expresses the class of the  pull back by the Abel-Jacobi map of the
theta
divisor on the Jacobian in terms of the classes  $x$ and $\Delta $.
That is,
$u^{*} \theta =
(d+g-1)x - \frac{\Delta }{2}$.  \\
\\
\noindent {\bf {Construction of  (normalized)  ${\cal L}_d$:}}
Consider the
Abel-Jacobi map $u :\, \univ ^{(d)} \lra \jac ^d $.  The bundle
$\frac{d+g-1}{s_d} {\cal L}_{\omega } -k_d \frac{\Delta }{2}$ is
trivial on the
fibers of $u$: Indeed, by Remark \ref{rem21}, we have  ${\cal
L}_{\omega
}|_{u^{-1}(L)} \cong {\cal O}(2g-2)$ and, by MacDonnald's formula, we
have
$\frac{\Delta}{2}|_{u^{-1}(L)} \cong {\cal O}(d+g-1)$.  We thus get
by the see
saw principle, see \cite{Mu1}, that the above bundle descents to a
bundle on
$\jac ^d$  and this is exactly the bundle  ${\cal L}_d$.  That ${\cal
L}_d$ has
class  $k_d \theta$, it is again an  application of the  MacDonnald's
formula.\\
\\
\noindent {\bf {Construction of   ${\cal P}_d$: }}  Consider the
diagram
\begin{equation}
\begin{diagram}[aaaaaaaaa]
 \node{\univ ^{(d)} \times_{\mod} \univ } \arrow{e,t}{\tilde {u}}
\arrow{s,r}{
p_1 }
 \node{\jac ^d \times_{\mod } \univ }  \arrow{s,r} {{\nu }} \\
\node{\univ ^{(d)} } \arrow{e,t}{u  } \node{\jac ^d}
\end{diagram}
\label{diag3}
\end{equation}

On $\univ ^{(d)} \times_{\mod } \univ $ consider the line bundle
$p_1^{*}(m{\cal L}_{\omega } -n \frac{\Delta }{2} )+ s_d {\cal D}$,
where $n,\,
m$ are integers satisfying $n(d+g-1)-m(2g-2)=s_d$.  One can see again
that this
is trivial on the fibers of $\tilde {u}$ (note that ${\cal
D}|_{u^{-1}(L)
\times \{ p \} }\cong {\cal O}(1)$) and that  its restriction on  $\{
D \}
\times C$ is  isomorphic  to $D^{\otimes s_d}$. The bundle ${\cal
P}_d $ is the
descent of that one on $\jac ^d \times_{\mod } \univ $.  Note that
different
choices of $n$, give rise to line bundles ${\cal P}_d$ which differ
by the pull
back of a line bundle from $\jac ^d$. In the following we can fix one
such $n$,
for example  choose $n$ to be the smallest positive integer such that
the
number $\frac{n(d+g-1) - s_d}{2(g-1)}$ is an
 integer, or equivalently,  such that the number $\frac{nd-s_d}{g-1}
+ n$  is
an  even integer.

We continue with some lemmas about intersections. In the following we
denote a
line bundle and its first Chern class in the Chow group by the same
symbol.

\begin{lem}
If $p_1 : \, \univ ^{(d)} \times_{\mod} \univ \lra \univ ^{(d)}$
denotes the
first projection, then
$$ p_{1*}{\cal D}^2 = -{\cal L}_{\omega } + \Delta  . $$
\label{D}
\end{lem}
\vskip-.2in
 {\sc {Proof:}} Consider the diagram

\begin{equation}
\begin{diagram}[aaaaaaaa]
 \node{\univ ^{\times d}} \arrow{e,t}{\delta _{i,d+1}}
\arrow{se,t}{id}
 \node{\univ ^{\times d} \times_{\mod } \univ} \arrow{e,t}{\tilde
{c}}
 \arrow{s,r}{ q_1 }
 \node{\univ ^{(d)} \times_{\mod } \univ }   \arrow{s,r}{p_1} \\
 \node{\univ } \node{\univ ^{\times d}}  \arrow{w,t}{\pi _i}
\arrow{e,t}{c }
  \node{\univ ^{(d)}}
\end{diagram}
\label{diag4}
\end{equation}
\noindent
The map $\delta_{i,d+1}$ is the $i,d+1$-diagonal embedding, and $c$
the
canonical map. The maps $q_1$  and $\pi _i$ are the projections.
The map $c$ is flat and its pull back defines an injection in the
intersection
rings. By the commutativity and Remark \ref{rem22}, it is enough to
show that
$q_{1*} \tilde{c}^{*}({\cal D} ^2)= -{\cal L}_K +2\Delta$, where we
denote also
 by $\Delta $ the sum of the big diagonals in $\univ ^{\times d}$. We
have
\[
\tilde{c}^{*}({\cal D} ^2)=\tilde{c}^{*}({\cal D})^2 = ( \sum
_{i=1}^{d}
\Delta _{i,d+1})^2 = \sum_{1 \leq i,j \leq d} \Delta _{i,d+1} \,
\Delta _{j,
d+1}.
\]
 There are two cases to consider
 \begin{eqnarray*}
 \mbox{if} \;\; & i\neq j & \;\; \mbox{then} \;\;\; q_{1*}( \Delta
_{i,d+1}
\Delta_{j, d+1})=\Delta _{i,j}\\
 \mbox{if} \;\; & i=j & \;\; \mbox{then} \;\;\; q_{1*} \Delta
_{i,d+1}^2 =
\delta _{i,d+1}^{*} \Delta _{i,d+1} = \pi _i^{*} \omega _{\pi
_i}^{-1}.
 \end{eqnarray*}
 Therefore
 \begin{eqnarray*}
 q_{1*} \tilde{c}^{*}({\cal D} ^2) & = & 2 \sum _{1\leq i<j \leq
d}\Delta
_{i,j} + \sum_{i=1}^d \omega _{\pi _i}^{-1} \\
    & = & 2\Delta - {\cal L}_K.
 \end{eqnarray*}
\begin{flushright}
   $\Box$
   \end{flushright}

 \begin{lem}
 Following the notation of diagram \ref{diag3}, we have the following
(where
again we keep the same notation for a line bundle and its first Chern
class):
\begin{enumerate}
\item  $\nu _{*}( {\cal P }_d ^2) = s_d^2  \,
\frac{nd-s_d}{g-1}\,{\cal L }_d
$,
 \item $\nu _{*} ({\cal P }_d \omega_{\nu }) = s_d n \, {\cal L }_d$,
\item $\nu _{*}(\omega _{\nu }^2) = 12 \, \psi ^*\lambda$,
\end{enumerate}
where $n$  is the integer defined in  the construction of ${\cal
P}_d$, see
Section 2.
\label{push forward}
\end{lem}
{\sc {Proof:}} Again, since $u$ is flat and its pull back defines an
injection
in the chow rings, it is enough to prove that $p_{1*} \tilde {u} ^{*}
({\cal
P}_d^2) = s_d^2 \,
\frac{nd-s_d}{g-1}\, u^*{\cal L}_d$. By the construction of ${\cal
P}_d$, we
have
\[
 \tilde {u} ^{*}
 ({\cal P}_d^2)   = \left[ s_d\,  {\cal D}+ p_1^{*}(m {\cal
L}_{\omega } - n
\frac{\Delta }{2}) \right] ^2= s_d^2  \, {\cal D}^2+ 2 s_d  \, {\cal
D} \,
p_1^{*} (m {\cal L}_{\omega } - n \frac{\Delta }{2}) + p_1^{*}(m
{\cal
L}_{\omega } - n \frac{\Delta }{2})^2.
 \]
  Thus,
  \begin{eqnarray*}
  p_{1*} \tilde{u}^{*}{\cal P}_d^2 & = & s_d^2 \,  p_{1*}{\cal D}^2 +
2 s_d\,
p_{1*} {\cal D}
  (m {\cal L}_{\omega }  - n \frac{\Delta }{2})
   =  s_d^2 \, (-{\cal L}_{\omega }+\Delta ) +
 2 d s_d  \,  (m {\cal L}_{\omega } - n \frac{\Delta }{2}) \\
  & = & (2 d m  s_d   -s_d ^2)
 {\cal L}_{\omega } + (2 s_d^2 - 2 d n s_d ) \frac{\Delta }{2}.
 \end{eqnarray*}
  The first coefficient can be written as
  \begin{eqnarray*} 2 d m s_d
 -s_d ^2 & = & \frac{d s_d }{g-1} (n(d+g-1)-s_d)-s_d^2 = \frac{d s_d
}{g-1}  \,
n (d+g-1)- \frac{d s_d ^2 }{g-1}- s_d ^2 \\
  & = &  \frac{d s_d }{g-1} \,  n (d+g-1)- \frac {s_d ^2}{g-1}(d+g-1)
=
   \frac{s_d }{g-1}(d+g-1)(dn-s_d).
  \end{eqnarray*}
   The second coefficient can be written as
   \[
   2 s_d^2 - 2 d n s_d  = -2s_d (dn -s_d).
   \]
    We thus have
    \begin{eqnarray*}
 p_{1*} \tilde {u}^*  {\cal P}_d^2  & = &  \frac{s_d}{g-1}
(d+g-1)(dn-s_d) \,
 {\cal L}_{\omega }-2s_d(dn-s_d)\frac{\Delta }{2}\\
  & =  & \frac{s_d^2}{g-1}(dn-s_d) \frac{d+g-1}{s_d} \,  {\cal
L}_{\omega } -
s_d ^2 \frac{dn-s_d}{g-1} k_d \, \frac{\Delta }{2} \\
  & = & s_d ^2 \frac{dn-s_d }{g-1} \,  (\frac{d+g-1}{s_d } \, {\cal
L}_{\omega
} -k_d \frac{\Delta }{2})=  s_d ^2 \frac{dn-s_d }{g-1}  \, u^{*}
{\cal L}_d.
 \end{eqnarray*}
  This proves the first part of the lemma.

 For the second part, consider the diagram
 \vskip-.3in
\begin{equation}
\begin{diagram}[aaa]
  \node{ } \node[2] { } \node[2]{\univ }  \\
 \node { }
 \node{\univ ^{(d-1)} \times_{\mod } \univ }
   \arrow[2]{e,t}{\delta } \arrow{sw,t}{\pi _2} \arrow{se,t}
   {\sigma }
 \node[2]{\univ ^{(d)} \times_{\mod } \univ }  \arrow{ne,t}{p_2}
\arrow
{sw,t}{p_1} \arrow[2]{e,t}{\tilde {u}} \node[2]{\jac ^d \times _{\mod
} \univ}
\arrow{sw,t}{\nu }
   \\
 \node{ \univ }
 \node[2]{\univ ^{(d)} }\arrow[2]{e,t}{u}
 \node[2]{\jac ^d}
\end{diagram}
\label{diag5}
\end{equation}
\noindent
The map $\sigma $ in the diagram is the addition map. We have
\[
\tilde {u}^{*} ({\cal P}_d \omega _{\nu }) = (p_1^{*}(m{\cal
L}_{\omega } - n
\frac{\Delta }{2}) + s_d {\cal D}) \omega _{p_1 } = p_1^{*}(m{\cal
L}_{\omega }
- n \frac{\Delta }{2}) p_2^{*} \omega _{\pi } + s_d  \, {\cal D}
p_2^{*} \omega
_{\pi }.
\]
By  applying $p_{1*}$ and using the definition of ${\cal D}$, we get
 \[
 p_{1*} \tilde{u}^{*} ({\cal P}_d \omega _{\nu })= (2g-2) (m{\cal
L}_{\omega }
- n \frac{\Delta }{2}) + s_d \, \sigma _{*} \delta ^{*} p_2^{*}
\omega _{\pi }.
 \]
   Since $p_2 \circ \delta = \pi _2$, we have that
  $ \sigma _{*} \delta ^{*} p_2^{*} \omega _{\pi }= \sigma _* \pi
_{2}^* \omega
_{\pi }$.
  By the definition of ${\cal L}_{\omega}$  and by some diagram
chasing
   it is easy to see that
  $ \sigma _{*}\pi _2 ^{*} \omega _{\pi }= {\cal L}_{\omega }$.  We
thus get
  \begin{eqnarray*}
  p_{1*} \tilde{u} ^{*} ({\cal P}_d \omega _{\nu }) & = & ((2g-2)m
+s_d ){\cal
L}_{\omega } - n (2g-2)\frac{\Delta }{2} \\
   & = &  n(d+g-1) {\cal L}_{\omega } - n (2g-2)\frac{\Delta }{2} = n
s_d  \,
u^{*}{\cal L}_d.
   \end{eqnarray*}
  This proves the second part.

 The third part is an immediate consequence of the  fact $\pi _{*}
(\omega
_{\pi }^2) = 12 \lambda $, see \cite{Mu2}.
\begin{flushright}
   $\Box$
   \end{flushright}

  \section{The Hilbert scheme}
Let $C$ be  a  smooth curve of genus $g \geq 4$ and $d>2g-2 $, $r
\leq d-g $
given  integers. We give now estimates for the dimension of the
complement of
$\hilb $ in ${\frak P}_d$. Let $L$ denote a line bundle of degree $d$
on $C$.
The complement of $\hilb $ in the fiber
${\frak P}_d^L = \pr ({\Bbb C}^{r+1} \otimes  H^0(C,L))$  of
${\frak P}_d \lra \jac ^d$ over  $L$, consists of two (maybe
overlapping) loci
${\frak U}^L_{\rm deg}$ and ${\frak U}^L_{\rm nemb}$.   The first
contains
those points which define degenerate maps and the second those for
which the
corresponding map is not an embedding of degree $d$.  We have the
following
lemmas:

\begin{lem}
Following the above notation, we have that ${\rm codim}_{{\frak
P}_d^L}
{\frak U}^L_{\rm deg } = d-g-r+1 $.
 In particular, if $3 \leq r < d-g $, then ${\rm codim}_{{\frak
P}_d^L}
{\frak U}^L_{\rm deg }  \geq 2$ and if $r= d-g $, then ${\rm
codim}_{{\frak
P}_d^L}
{\frak U}^L_{\rm deg } =1$. In the later case, ${\frak U}^L_{\rm deg
}$ is an
irreducible divisor of degree $d-g+1$ in the projective space
${\frak P}_d^L$.
\label{complement1}
\end{lem}
{\sc {Proof:}} The above locus ${\frak U}^L_{\rm deg} $ corresponds
to those
tuples
  $<s_1, \ldots , s_{r+1}> \in
  {\Bbb C}^{r+1} \otimes  H^0(C,L)$ which span a space of dimension
$\leq r$.
  Since  dim$H^0(C,L)=d-g+1$, this proves the Lemma.

\begin{flushright}
   $\Box$
   \end{flushright}

\begin{rem}
{\rm  For $r=d-g$, the assumption that $d \geq  2g+1$, implies that $
{\frak
U}^L_{\rm nemb} \subset {\frak U}^L_{\rm deg} $  with  codimension
$\geq 1$.
 }
\label{cased-g}
\end{rem}

\begin{lem}
For $ 4 \leq r < d-g $, we have that $ {\rm codim}_{{\frak P}_d^L}
 {\frak U}^L_{\rm  nemb} \geq  2 $.
For $3=r < d-g$, the locus $ {\frak U}^L_{\rm nemb}  $  is an
irreducible
divisor of degree $2(d-1)(d-2)-4g$  in the projective space ${\frak
P}_d^L  $.
\label{grass}
\end{lem}
{\sc {Proof:}}  The space ${\frak P}_d^L \setminus {\frak U}^L_{\rm
deg}$
 maps in a natural way to the Grassmanian ${\bf Gr}(r+1,H^0(C,L))$
which
parametrizes linear series $g^r_d$ of $L$ on $C$.  Thus, for $r <
d-g$,  there
is a rational map $\alpha :  {\frak P}_d^L  \lra {\bf
Gr}(r+1,H^0(C,L))$  which
is not defined in a codimension  $\geq 2$ locus, see Lemma
\ref{complement1}.
The fiber of $\alpha $ is isomorphic to
  ${\bf PGL}(r+1)$.  The locus ${\frak U}^L_{\rm nemb} \setminus
  {\frak U}^L_{\rm deg} $ is the pull back of the correspondent locus
in
  ${\bf Gr}(r+1,H^0(C,L))$.

 It is enough to prove the Lemma on the ``level'' of  ${\bf
Gr}(r+1,H^0(C,L))$.
  Consider
  the curve $C$ embedded in $\pr (H^0(C,L)^{\vee})$ by the complete
linear
system of $L$.
  Maps defined by the $g^r_d$'s, correspond to projections from
$(d-g-r-1)$-dimensional
  projective planes in the above space.
 Maps which  are not embeddings of degree $d$, correspond
 to projections from those planes which  intersect the secant
variety of $C$.
The later is a locus in the dual of the Grassmanian  ${\bf
Gr}(r+1,H^0(C,L))$
of codimension  equal to
 $ r-2 $ (the dimension of the secant variety is 3).  Thus,  if
$r\geq 4$, then
the codimension  of  $ {\frak U}^L_{\rm nemb}  $  is $\geq 2$.

 We turn now to the case $r=3$.  Observe first that the pull back by
$\alpha $
of the generator ${\cal O}_{\bf Gr}(1)$ of the Picard group of the
Grassmanian
to
  ${\frak P}_d^L$ is isomorphic to ${\cal O}_{{\frak P}_d^L}(r+1)$.
By the
above
 discussion,  one can see that ${\frak U}^L_{\rm nemb}$   is  an
irreducible
divisor;
  the formula
 for its degree  is a consequence of the previous observation and  of
the fact
that the degree
  of the secant variety is $\frac{(d-1)(d-2)}{2}-g$.
\begin{flushright}
   $\Box$
   \end{flushright}
\bigskip
 We have the following lemmas. The proof of the first is an immediate
consequence of Lemmas \ref{complement1} and  \ref{grass} above.

 \begin{lem}
For $d \geq 2g+1$ and $3<r<d-g$,  we have an isomorphism of Picard
groups
$\mbox{Pic} \, \hilb \cong \mbox{Pic} \, {\frak P}_d$.
 \label{equal picard}
 \end{lem}

 \begin{lem}
 The bundle ${\frak P}_d \lra \jac ^d$ admits a line bundle whose
restriction
to a fiber is isomorphic to ${\cal O}(s_d)$. The number $s_d $ is
minimum with
this property.
 \label{number}
 \end{lem}
{\sc {Proof:}} Lets consider first the same question for the bundle
$u:\, \symm
\lra \jac ^d $. Let ${\cal L}$ be a line bundle on $\symm $ whose
restriction
on the fiber $\pr H^0(C,L)$ is isomorphic to ${\cal O}(t)$.  By
\cite{K1}, pg.
844 , the class of its restriction
  is of  the form $(2g-2)m x - n \frac{\Delta }{2}$, where $m, \, n$
are
integers.  Since the restrictions of $x$ and $\frac{\Delta}{2} $ on a
fiber of
the Abel-Jacobi map  have classes $c_1{\cal O}(1)$ and
  $c_1{\cal O}(d+g-1)$ respectively,   we conclude  that $s_d|t$.
One can
see that the same is true for the bundle ${\frak P}_d$. The proof is
similar to
that of Lemma 4, in \cite{K2}.

On the other hand, one  can construct  a  line bundle ${\cal R}$ on
$\hilb $
  whose restriction on the fibers is isomorphic to ${\cal O}(s_d)$ as
follows.
We turn  back to diagram \ref{diag1}. The pull back by the map $\eta
$ of
${\cal O}(1)$  to $\unhilb $,  is the tautological bundle ${\cal F}$.
By the
see-saw principle we have that there exists a bundle ${\cal R}$ on
$\hilb $
such that
  \[
  q^{*} {\cal R} \cong {\cal F}^{\otimes s_d} \otimes {\alpha }^{*}
{\cal P}_d
^{-1}.
  \]
It is easy to see that the restriction of ${\cal R}$ to the fibers of
the map
$\phi $ is ${\cal O}(s_d)$.
\begin{flushright}
   $\Box$
   \end{flushright}
\bigskip

  The above constructed line bundle ${\cal R}$  is the generator of
the
relative Picard group  of $\hilb$ over $ \jac ^d$. To summarize,  we
have that
the Picard group of
  $\hilb $, when $d \geq 2g+1 $ and $3<r<d-g $,  is freely generated
over the
integers
  by the  line bundles $\psi ^{*} \phi ^{*} \lambda, \; \phi ^{*}
{\cal L}_d$
and
  ${\cal R}$.

\section{The classes of the generators}
In this section we express the classes of the above generators of
Pic$\hilb $,
in terms of the naturally defined classes $A, \; B, \; C $ of Section
1. For
the notation, see diagram 1.

\begin{lem}
If   $A = q_{*}({\cal F}^2), \;\;  B  = q_{*} ({\cal F} \omega _q)$
and $ C
=q_{*}( \omega _q^2)$, we have
\begin{enumerate}
 \item $ A = \frac {nd-s_d}{g-1} \, \phi  ^{*} {\cal L}_d + 2
\frac{d}{s_d }
{\cal R} $,
 \item $ B= n \, \phi ^{*}  {\cal L}_d + \frac{2g-2}{s_d}\,  {\cal R}
$,
 \item $ C  =  12 \,   \phi ^{*} \psi ^{*} \lambda $,
\end{enumerate}
 where $n$ is the integer defined in the construction of ${\cal P}_d$
in
Section 2.
\label{classes}
\end{lem}
\begin{rem}
{\rm  Note that the coefficient matrix has determinant  of absolute
value 24.}
\end{rem}
\noindent
{\sc {Proof:}}  For the first one:
\begin{eqnarray*}
  s_d^2 \, {\cal F}^2  & = & (q^{*}{\cal R} + \alpha ^{*} {\cal
P}_d)^2 \\
              &  = &  q^{*}{\cal R}^2 + 2 \,  q^{*}{\cal R} \,
\alpha ^{*}
{\cal P}_d +
	       \alpha ^{*} {\cal P}_d^2.
\end{eqnarray*}
Therefore by applying  $q_{*}$, we get
\begin{eqnarray*}
s_d^2 \, q_{*} {\cal F}^2  & = & 2 \, q_{*} \alpha ^{*} {\cal P}_d \,
{\cal R}
+ q_{*} \alpha ^{*}
  {\cal P}_d^2  \\
    & = & 2 d s_d  \, {\cal R} + \phi ^{*} \nu _{*} {\cal P}_d^2   \\
    &  \stackrel {\mbox{Lem \ref {push forward}}}{=} &  2  d s_d  \,
{\cal R} +
s_d^2 \,  \frac{nd-s_d}{g-1} \phi ^{*} {\cal L}_d.
 \end{eqnarray*}
 and this proves the first part.

 For the second part: By the definition of  ${\cal R}$ and by
multiplying by
 $\omega _q$ we have
\[
{\cal F} \omega _q = \frac{1}{s_d} \, q^{*}{\cal R} \, \omega _q   +
\frac{1}{s_d} \, \alpha  ^{*} {\cal P}_d  \, \omega _q,
\]
and so,
\begin{eqnarray*}
q_{*}({\cal F} \omega _q)  & = & \frac{1}{s_d} \, {\cal R} \, q_{*}
\omega _q
  + \frac{1}{s_d} \,  q_{*}  \alpha  ^{*}( {\cal P}_d \omega _{\nu })
\\
 & = & \frac{2g-2}{s_d} \, {\cal R} + \frac{1}{s_d} \, \phi ^{*} \nu
_{*}
({\cal P}_d \omega _{\nu })  \\
 &  \stackrel {\mbox{Lem \ref {push forward}}}{=} &
\frac{2g-2}{s_d}\,  {\cal
R} + n \phi ^{*}\,  {\cal L}_d.
\end{eqnarray*}

The third part is an immediate consequence  of Lemma \ref{push
forward}.
\begin{flushright}
   $\Box$
   \end{flushright}
\bigskip

The above Lemma  \ref{classes} gives the expression of the generators
in terms
of the classes $A,\, B, \, C $. We thus have

  \begin{theor}
 For $d \geq 2g+1 $ and $3 < r < d-g$, the Picard group of $\hilb $
is  freely
generated over $\Bbb Z$ by the three line bundles $\psi ^{*} \phi
^{*} \lambda,
\; \phi ^{*} {\cal L}_d$ and ${\cal R}$. Their classes can be
expressed as
 \begin{enumerate}
 \item  $ {\cal R}=  \frac {1}{2} n \,  A-  \frac {1}{2}\frac{nd -
s_d}{g-1}\,
B ,   $
 \item   $  \phi ^{*} {\cal L}_d =   - \frac {g-1}{s_d} \, A +
\frac{d}{s_d} \,
B , $
  \item   $  \phi ^{*} \psi ^{*} \lambda = \frac{1}{12}  \,  C   ,  $
 \end{enumerate}
 where $n$ is the integer defined in the construction of ${\cal
P}_d$, see
Section 2.
  \label{main}
\end{theor}
\begin{rem}
{\rm Note that in the above theorem the numbers $\frac{g-1}{s_d}$ and
  $  \frac{d}{s_d}  $ are either integers or half integers.  The
number $
\frac{nd - s_d}{g-1}$ is an integer  by its definition. }
\end{rem}
\begin{rem}
{\rm  The cases $r=3 < d-g $ and $r=d-g$ have to be treated
separately: For
$r=3 < d-g $
   the Picard group of ${\frak H}_{d,g,3}$  fits in an exact sequence
 \[
 0 \lra {\Bbb Z}  \stackrel {\mu }{ \lra}  {\Bbb Z}^{\oplus 3}
\stackrel {p
}{ \lra} \mbox{Pic} \, {\frak H}_{d,3,g}  \lra  0 ,
 \]
where the map $\mu $ is the multiplication by  $2(d-1)(d-2) -4g $ in
the first
factor.
 The case is similar for $r=d-g$. Here $\mu $ is the multiplication
by $d-g+1$.
We leave
  to the reader to find in this case the analogue of Theorem
\ref{main}.}
\end{rem}

\begin{rem}
{\rm   One can pursue the above discussion further and calculate the
Picard
groups over the integers of the Severi varieties and the Hurwitz
schemes in the
case when the degree $d$ is big with respect to the genus $g$.    For
the
Severi varieties, the results of Diaz-Harris, see \cite{DH1},
\cite{DH2}, imply
that the union of the Severi variety $V_{d,g}$ with the three
irreducible,
independent divisors $CU, \; TN, \; TR$ of ${\frak P}_d$ which are
defined in
the above mentioned papers, has complement of codimension $\geq 2$ in
${\frak
P}_d$. Since Pic${\frak P}_d$ is isomorphic to ${\Bbb Z}^{\oplus 3}$,
this
implies that Pic$V_{d,g}$ is torsion. Furthermore, by using the
expressions of
$CU, \; TN, \; TR$ in terms of $A,B,C$ given  in \cite{DH1} and the
above Lemma
\ref{classes}, one can calculate that torsion group. A similar result
should be
obtained for the Hurwitz scheme by using the analogues results of
Mockizuki,
see \cite{Mo}.  In the case of the Severi varieties, the calculations
of the
author lead to rather messy formulas.    }
\end{rem}

\end{document}